\documentclass[final,leqno]{siamltex}
\usepackage{amsmath}
\usepackage{upref} 
\usepackage{amssymb}
\usepackage{breakcites}
\usepackage{color}
\usepackage{array}

\usepackage{graphicx}
\usepackage{url}
\usepackage{placeins} 
\usepackage{fancyhdr,substr}  
\usepackage{gitinfo}

\definecolor{pred}{rgb}{0.8,0.1,0}
\definecolor{pgreen}{rgb}{0.3,0.7,0}
\definecolor{pblue}{rgb}{0.1,0,1.0}
\definecolor{hotpink}{rgb}{0.9,0,0.5}

\usepackage[colorlinks,urlcolor=pred,citecolor=hotpink,linkcolor=pblue]{hyperref}

\makeatletter
\def\@cite#1#2{[{#1\if@tempswa ,~#2\fi}]}
\makeatother
\title{A New Algorithm for Computing the Actions of Trigonometric
and Hyperbolic Matrix Functions%
       \thanks{Version of \today}
}

\author{Awad H. Al-Mohy%
           \thanks{%
                   Department of Mathematics, King Khalid University, Abha,
                   Saudi Arabia
                   (ahalmohy@kku.edu.sa,
                    http://www.maths.manchester.ac.uk/\string~almohy).
        }
}

\makeatletter
\def\argmin{\mathop{\operator@font argmin}}
\makeatother

\def\mv{mv}

\def\speed{$t_{\mathrm{ratio}}$}

\def\shift{\mbox{\upshape{shift}}}
\def\parameters{\mbox{parameters}}

\def\dfrac#1#2{{\displaystyle{#1\over#2}}}
\def\matvec{matrix--vector}

\makeatletter
\def\mycases#1{\left\{\,\vcenter{\normalbaselines\m@th
    \ialign{$##\hfil$&\quad{##}\hfil\crcr#1\crcr}}\right.}
\makeatother

\def\Cmstar{C_{m_*}}
\def\pmax{p_{\max}}

\def\undoin{\t{undoin}}
\def\undout{\t{undout}}

\newtheorem{algorithm}[theorem]{Algorithm}
\newtheorem{fragment}[theorem]{Code Fragment}

\def\half{\frac{1}{2}}

\newcommand{\parens}[1]{\rom{(}#1\rom{)}}

\def\At{\widetilde{A}}

\def\xor{\mathrm{~xor~}}
\def\C{\mathbb{C}}

\def\and{\mathop{\mathrm{and}}}
\def\sinch{\mathop{\mathrm{sinch}}}

\def\half{\frac{1}{2}}

\def\mmax{m_{\max}}
\def\pmax{p_{\max}}

\def\mystrut#1{\rule{0cm}{#1}}  

\def\sqrtA{A^{1/2}}

\def\Pant{Pad\accent"13e approximant}

\mathcode`@="8000 
{\catcode`\@=\active\gdef@{\mkern1mu}}

\newcounter{Muni}
\renewcommand{\theMuni}{\alph{Muni}}

\def\rom#1{{\upshape#1}}

\def\rt{\widetilde{r}}

\def\alg{algorithm}
\def\Alg{Algorithm}

\def\t#1{\texttt{\upshape #1}}

\def\fpa{floating point arithmetic}

\def\Myfig#1{\begin{center}\includegraphics[width=10.5cm]{#1}\end{center}}

\newcounter{mylineno}
\makeatletter
\let\oldtabcr\@tabcr

\def\mynewline{\refstepcounter{mylineno}%
                \llap{\footnotesize\arabic{mylineno}\hspace{5pt}}%
               }

\gdef\@tabcr{\@stopline \@ifstar{\penalty%
            \@M \@xtabcr}\@xtabcr\mynewline}

\newenvironment{code}{%
                         \mathcode`\:="603A  
                         \def\colon{\mathchar"303A}
                         \setcounter{mylineno}{0}
                         \par
                         \upshape
                         \begin{list} 
                            {} {\leftmargin = 1.8cm}
                         \item[]
                         \begin{tabbing}

                            \hspace*{.3in} \= \hspace*{.3in} \=
                            \hspace*{.3in} \= \hspace*{.3in} \= \kill
                            \mynewline
                        }{\end{tabbing}\end{list}}
\makeatother
\def\a{\alpha}

\def\l{\lambda}
\def\sg{\sigma}
\def\th{\theta}

\def\resp{respectively}

\def\xhat{\widehat{x}}

\def\tol{\mathrm{tol}}
\def\cond{\mathrm{cond}}

\def\sinc{\mathrm{sinc}}

\def\trace{\mathrm{trace}}
\def\diag{\mathrm{diag}}
\def\d{\mathrm{d}}
\def\i{\mathrm{i}}

\def\Re{\mathrm{Re}}
\def\Im{\mathrm{Im}}
\def\norm#1{\|#1\|}
\def\normt#1{\|#1\|_2}

\def\normi#1{\|#1\|_1}
\def\normo#1{\|#1\|_{\infty}}

\def\nbyn{n \times n}

\def\C{\mathbb{C}}
\def\R{\mathbb{R}}

\let\oldref\ref
\def\ref#1{{\normalfont\oldref{#1}}}
\def\eqref#1{{\normalfont(\oldref{#1})}}

\makeatletter
\def\mymatrix#1{\null\,\vcenter{\normalbaselines\m@th
    \ialign{\hfil$##$\hfil&&\quad\hfil$##$\hfil\crcr
      \mathstrut\crcr\noalign{\kern-\baselineskip}
      #1\crcr\mathstrut\crcr\noalign{\kern-\baselineskip}}}\,}
\makeatother

\newcounter{example}
\def\Experiment{Experiment}

\newenvironment{example}%
               {\begin{list}{\indent\emph{\Experiment} {\upshape\arabic{example}.}}%
                {\usecounter{example}
                \setlength{\leftmargin}{\rightmargin}
                \setlength{\labelwidth}{\leftmargin}
                \addtolength{\labelwidth}{-\labelsep}
                \setlength{\topsep}{0in}
                \setlength{\itemsep}{0pt}
                \setlength{\listparindent}{\parindent}%
               }}%
               {\end{list}}

\mathchardef\Gamma="7100 \mathchardef\Delta="7101
\mathchardef\Theta="7102 \mathchardef\Lambda="7103
\mathchardef\Xi="7104 \mathchardef\Pi="7105 \mathchardef\Sigma="7106
\mathchardef\Upsilon="7107 \mathchardef\Phi="7108
\mathchardef\Psi="7109 \mathchardef\Omega="710A

\begin{document}
\maketitle

\begin{abstract}
A new algorithm is derived for computing the actions $f(tA)B$ and
$f(tA^{1/2})B$, where $f$ is cosine, sinc, sine, hyperbolic cosine,
hyperbolic sinc, or hyperbolic sine function. $A$ is
an $n\times n$ matrix and $B$ is $n\times n_0$ with $n_0 \ll n$.
$A^{1/2}$ denotes any matrix square root of $A$ and it is never required
to be computed.
The algorithm offers six independent output options given $t$, $A$, $B$, and
a tolerance. For each option, actions of a pair of trigonometric or
hyperbolic matrix functions are simultaneously computed.
The algorithm scales the matrix $A$ down by a positive integer $s$,
approximates $f(s^{-1}tA)B$ by a truncated Taylor series, and finally
uses the recurrences of the Chebyshev polynomials of the first and
second kind to recover $f(tA)B$. The selection of the scaling parameter
and the degree of Taylor polynomial are
based on a forward error analysis and a sequence of the form $\|A^k\|^{1/k}$
in such a way the overall computational cost of the algorithm is optimized.
Shifting is used where applicable as a preprocessing step to reduce
the scaling parameter.
The algorithm works for any matrix $A$ and
its computational cost is dominated by the formation of products of
$A$ with $n\times n_0$ matrices that could take advantage of the
implementation of level-3 BLAS.
Our numerical experiments show that the new algorithm behaves
in a forward stable fashion and in most problems outperforms
the existing algorithms in terms of CPU time, computational cost,
and accuracy.

\end{abstract}




\section{Introduction}

The matrix cosine and sine functions appear in the solution
of the system of second order differential equations
\begin{equation}\label{ode}
  \frac{\d^2y}{\d t^2} + Ay = g(y(t)),\quad y(0)=y_0,\quad y'(0)=y'_0.
\end{equation}
The exact solution of this system and its derivative
is given by the variation of the constants formula \cite{grho06,wu15}

\begin{eqnarray}
\label{sol.ode}
 y(t) = \cos(t\sqrtA)y_0 &+& t\,\sinc(t\sqrtA)y'_0\\
&+&\int_{0}^{t}(t-\tau)\,\sinc((t-\tau)\sqrtA)g(y(\tau))\d \tau,
\nonumber
\end{eqnarray}
\begin{eqnarray}
\label{sol.ode.y'}
y'(t) = -\sqrtA\sin(t\sqrtA)y_0 &+& \cos(t\sqrtA)y'_0\\
&+&\int_{0}^{t}(t-\tau)\cos((t-\tau)\sqrtA)g(y(\tau))\d \tau,
\nonumber
\end{eqnarray}
where $\sqrtA$ denotes any matrix square root of $A$
and $\sinc :\C^{\nbyn}\to\C^{\nbyn}$ is defined as
\begin{equation}
\label{sinc.series}
  \sinc X = \sum_{k=0}^{\infty}\frac{(-1)^kX^{2k}}{(2k+1)!}.
\end{equation}
The matrix function $\sinc$ clearly satisfies the relation
 $X\sinc X = \sin X$.
The first term of \eqref{sol.ode.y'} can be rewritten using the quality
$$
\sqrtA\sin(t\sqrtA) = tA\,\sinc(t\sqrtA).
$$
This is important to clear any ambiguity that a square
root of $A$ is needed.
We will see below how the actions of $\cos(t\sqrtA)$ and
$\sinc(t\sqrtA)$ can be simultaneously computed without
explicitly computing $\sqrtA$ whereas it is impossible to evaluate
the action of $\sin(t\sqrtA)$ without forming $\sqrtA$ explicitly
because $\sin$ is an odd function.

The variation of the constants formula forms the basis of numerical
schemes to solve the problem. For instance, at time $t_n=nh$, $y(t_n)$
and $y'(t_n)$ can be numerically approximated by $y_n$ and $y'_n$, \resp,
via the trigonometric scheme
\begin{eqnarray}
\label{yn}
\qquad\quad
   y_{n+1} &=& \cos(h\sqrtA)y_n+h\,\sinc(h\sqrtA)y'_n
           +\frac{h^2}{2} \sinc(h\sqrtA)\widehat{g}(y_n),\\
\label{y'n}
   y'_{n+1} &=& -hA\,\sinc(h\sqrtA)y_n+\cos(h\sqrtA)y'_n
           +\frac{h}{2} \cos(h\sqrtA)\widehat{g}(y_n)
           +\frac{h}{2}  \widehat{g}(y_{n+1}),
\end{eqnarray}
where $\widehat{g}(y) = \psi(h\sqrtA)g(\phi(h\sqrtA)y)$
provided that $\psi$ and $\phi$ are suitably chosen continuous filter functions;
see \cite[sect.~2]{glmrs17}, \cite[sect.~2]{grho06}, or
\cite[sect.~XIII.2.2]{hlw06}.
Many filter functions are proposed in literature and most of them
involve several actions of $\sinc(h\sqrtA)$ to evaluate
$\widehat{g}(y)$. For example Hairer and Lubich \cite{halu00} chose $\psi=\sinc$
and $\phi=1$ while Griman and Hochbruck proposed $\psi=\sinc^2$
and $\phi=\sinc$ \cite{grho06}.

The system \eqref{ode} arises from semidiscretization of some
second order PDE's by finite difference or finite elements methods
\cite{serb79}.
The hyperbolic matrix functions: $\cosh A$, $\sinh A$, and $\sinch A$,
where $\sinch A = \sinc(\i A)$, arise in the solution of coupled
hyperbolic systems of PDE's \cite{jnpc03}. They also have an application
in communicability analysis in complex networks \cite{ehh08}.
The matrix $A$ is
usually large and sparse, so finding methods to compute
the action of these matrix functions on vectors are
so crucial to reduce computational cost.

The computation of
the action of the matrix exponential has received significant research
attention; see \cite{alhi11} and the references therein. However it is not
the case for trigonometric and hyperbolic matrix functions. A possible
reason is that the second order system \eqref{ode} can be presented in
a block form of a first order system of ODE's and the matrix
exponential is used to solve the problem as in \eqref{blk.exp} below.
Grimm and Hochbruck \cite{grho08} proposed the use of
a rational Krylov subspace method instead of the standard one for certain
problems to compute $\cos(t\sqrtA)b$ and $\sinc(t\sqrtA)b$. Recently,
Higham and Kandolf \cite{hika17} derived an \alg\ to compute
the action of trigonometric and hyperbolic matrix functions.
They adapted the existing \alg\ of Al-Mohy and Higham \cite{alhi11}, \t{expmv},
for computing the action of the matrix exponential so that the evaluation
of $\cos(A)B$ and $\sin(A)B$ (or $\cosh(A)B$ and $\sinh(A)B$)
requires the action of $e^A$ on the matrix
$[B,B]/2\in\C^{n\times 2n_0}$.

The calculation of $\cos A$ and $\sin A$ for dense $A$ of medium size
is will-studied.
Serbin and Blalock \cite{sebl80} proposed an algorithm for
$\cos A$. It begins by approximating $\cos(2^{-s}A)$ by
a Taylor or \Pant, where $s$ is a nonnegative integer, and then
applies the double angle formula
$\cos(2A)=2\cos^2(A)-I$ on the approximant $s$ times to recover the
original matrix cosine.
An \alg\ by Higham and Smith \cite{hism03} uses the [8/8] \Pant\
with the aid of a forward error analysis to specify the scaling
parameter $s$. Hargreaves and Higham \cite{hahi05c} develop an
algorithm with a variable choice of the degree of \Pant s based
on forward error bounds in such a way the computational cost is
minimized. They also derive an algorithm that computes $\cos A$
and $\sin A$ simultaneously. Recently, Al-Mohy et al.
\cite{ahr15} derive new backward stable \alg s for computing
$\cos A$ and $\sin A$ separably or simultaneously using \Pant s and
rational approximations obtained from \Pant s to the
exponential function. They use triple angle formula to have an
independent \alg\ for $\sin A$.
In spite of the fact that the \alg s based on the double and
triple angle formulas for computing $\cos A$ and $\sin A$, \resp,
prove great success, it doesn't seem that these formulas
can be adapted to compute the action of these matrix functions.

In this paper we derive a new \alg\ for computing the action
of the trigonometric and hyperbolic matrix functions of the form $f(tA)B$
and $f(t\sqrtA)B$ without computing $\sqrtA$. The form $f(t\sqrtA)B$
appears in the variation of constants formula \eqref{sol.ode}--\eqref{sol.ode.y'}.
In contrast, the \alg\ of Higham and Kandolf
cannot compute $f(t\sqrtA)B$ without explicitly computing $\sqrtA$,
which is completely impractical. Moreover, their \alg\ cannot immediately
return $\sinc(tA)B$ or $\sinch(tA)B$.

The paper is organized as follows. In section \ref{sec1} we exploit the
recurrences of the Chebyshev polynomials and explain how the actions of
trigonometric and hyperbolic matrix functions can be computed.
In section \ref{sec2} we present forward error analysis using
truncated Taylor series and computational cost analysis to determine
optimal scaling parameters and degrees of Taylor polynomials for various
tolerances.
Preprocessing by shifting and termination criterion are discussed in
section \ref{sec.prep}.
We write our \alg\ in section \ref{sec3} and then give numerical
experiments in section \ref{sec4}.
Finally we draw some concluding remarks in section \ref{sec5}.
\section{Computing the actions $\boldsymbol{f(tA)B}$ and
$\boldsymbol{f(t\sqrtA)B}$}\label{sec1}
In this section we exploit trigonometric formulas and derive
recurrences to computing the action of the matrix functions
$\cos X$, $\sinc X$, $\sin X$, $\cosh X$, $\sinch X$, and $\sinh X$
on a thin matrix $B$. For an integer $k$ we have
\begin{equation}\label{s}
  \cos(kX) + \cos((k-2)X) = 2\cos(X)\cos((k-1)X).
\end{equation}
Let $T_k(X,B) = \cos(kX)B$ and simply denote it by $T_k$, where $k\ge0$.
Thus we obtain the three term recurrence
\begin{eqnarray}\label{ss}
 T_k + T_{k-2} = 2\cos(X)T_{k-1}= 2\,T_1(X,T_{k-1}),\quad k\ge2.
\label{rec.cos}
\end{eqnarray}
Observe that \eqref{rec.cos} is the recurrence that generates
\emph{Chebyshev polynomials of the first kind for} $T_0=1$ and $T_1=x$
\cite{maha03}.
The heaviest computational work in the recurrence \eqref{rec.cos}
lies in $T_1(X,T_{k-1})$ for all $k\ge1$.
Let $r$ be a rational approximation to the cosine
function, which we assume to be good near the origin, and choose a positive
integer $s\ge1$ so that $\cos(s^{-1}A)$ is well-approximated by
$r(s^{-1}A)$. Thus
$$
T_1(s^{-1}A,T_{k-1})=\cos(s^{-1}A)T_{k-1}\approx r(s^{-1}A)T_{k-1}.
$$
 The recurrence \eqref{rec.cos} with $X=s^{-1}A$ yields
$$
T_s(s^{-1}A,B)=\cos(A)B.
$$
We choose for $r$ a truncated Taylor series
$$
r_m(x) = \sum_{j=0}^{m}\dfrac{(-1)^j x^{2j}}{(2j)!}
$$
and compute the matrix $V=r_m(s^{-1}A)B$ using consecutive
matrix products as shown by the next pseudocode.

\begin{fragment}\label{code.frag1}
\end{fragment}
\begin{code}
$V = B$\\
for \= $k=1:m$\\
\label{line.gamma}
\> $\beta = 2k$, $\gamma = 2k-1$ \\
\label{line.sqrt.frag}
\> $B = AB$\\
\> $B = (AB)\left(s^2\beta\gamma\right)^{-1}$\\
\> $V = V + (-1)^kB$\\
end
\end{code}
Similarly we approximate $\sinc\,x$ by truncating the Taylor series in
\eqref{sinc.series} as
$$
\rt_m(x) = \sum_{j=0}^{m}\dfrac{(-1)^j x^{2j}}{(2j+1)!}.
$$
The matrix $V:=\rt_m(s^{-1}A)B$ can be evaluated using Code Fragment
\ref{code.frag1} after replacing $\gamma$ in line \ref{line.gamma}
by  $\gamma = 2k+1$. To evaluate $r_m(s^{-1}\sqrtA)B$ or $\rt_m(s^{-1}\sqrtA)B$ , we only need to delete line
\ref{line.sqrt.frag} of Code Fragment \ref{code.frag1}.

Next, to compute $\sinc(A)B$ consider the three term recurrence
\begin{equation}\label{rec.sinc}
  U_k - U_{k-2} = 2T_k,\quad k\ge2, \quad U_0 = B,\quad
  U_1 =2T_1= 2\cos(X)B.
\end{equation}
It is the recurrence that yields \emph{the Chebyshev polynomials
of the second kind} \cite{maha03}.
By induction on $k$, it easy to verify that
\begin{equation}\label{sin.U}
\sin(X)U_{k-1} = \sin(kX)B.
\end{equation}
Assume for a temporarily fixed positive integer $q\ge2$ that
\eqref{sin.U} holds for all $k$ with $q\ge k\ge2$.
The inductive step follows from
\begin{eqnarray*}
  \sin(X)U_{q+1} &=& 2\sin(X)T_{q+1}+\sin(X)U_{q-1} \\
                 &=& 2\sin(X)\cos((q+1)X)B+ \sin(qX)B\\
                 &=& \bigl[\sin((q+2)X)-\sin(qX)\bigr]B + \sin(qX)B =\sin((q+2)X)B.
\end{eqnarray*}
Since \eqref{sin.U} holds for every $X$ we conclude that
\begin{equation}\label{sinc.U.formula}
 \sinc(X)U_{k-1} = k\,\sinc(kX)B.
\end{equation}
For $X=s^{-1}A$ the recurrences \eqref{rec.cos} and \eqref{rec.sinc} can
intertwine the computation of $T_s=\cos(A)B$ and $U_{s-1}$. The matrix
$\sinc(A)B$ can be recovered by computing the action
$\sinc(s^{-1}A)U_{s-1}\approx \rt_m(s^{-1}A)U_{s-1}$
that can be achieved by a single execution of Code Fragment
\ref{code.frag1} with
$V=U_{s-1}$ and $\gamma=2k+1$ in line \ref{line.gamma}. Observe
that the calculation of $U_{s-1}$ via \eqref{rec.sinc} involves
only $s-2$ additions of $n\times n_0$ matrices provided that $T_k$, $1\le k\le s$,
are already computed from \eqref{rec.cos}. Such operations
are negligible. However, we can save about the half of these
operations by observing that
\begin{equation}\label{ev.od.sum}
  \frac{1}{2}U_{s-1} =
         \mycases{T_1+T_3+T_5\cdots+T_{s-1},&  if $s$ is even, \cr
                  \frac{1}{2}T_0+T_2+T_4+\cdots+T_{s-1},&  if $s$ is odd,}
\end{equation}
which can be easily derived from \eqref{rec.sinc}.

Given the relations between trigonometric and hyperbolic functions,
we can replace $\cos$ in \eqref{rec.cos} and \eqref{rec.sinc} by $\cosh$
and replace $\sinc$ in \eqref{sinc.U.formula} by $\sinch$ so that
the recurrence relations return $\cosh(A)B$, $\sinch(A)B$, and
$\sinh(A)B$.

\section{Forward error analysis and computational cost analysis}\label{sec2}
We use the truncated Taylor series $r_m$ and $\rt_m$ to approximate the
$\cos$ and $\sinc$ functions, \resp. Given a matrix $A\in\C^{\nbyn}$ and
tolerance $\tol$, we need
to determine the positive integer $s$ so that
\begin{equation}\label{fwd.bound}
  \norm{\cos(s^{-1}A^\sg)-r_m(s^{-1}A^\sg)}\le \tol,
\end{equation}
where $\sg$ is either 1 or $1/2$. We have
$$
\cos(s^{-1}A^\sg)-r_m(s^{-1}A^\sg)=
\sum_{j=m+1}^{\infty}\dfrac{(-1)^j (s^{-1}A^\sg)^{2j}}{(2j)!}.
$$
By \cite[Thm.~4.2(b)]{alhi09a} and since the tail of
the Taylor series of the cosine is an even function, we obtain
\begin{eqnarray}
 \nonumber 
  \norm{\cos(s^{-1}A^\sg)-r_m(s^{-1}A^\sg)} &\le& \sum_{j=m+1}^{\infty}\dfrac{\a_p(s^{-1}A^\sg)^{2j}}{(2j)!} \\
                     &=& \cosh(\a_p(s^{-1}A^\sg))- \sum_{j=0}^{m}\dfrac{\a_p(s^{-1}A^\sg)^{2j}}{(2j)!}=:
\rho_m(\a_p(s^{-1}A^\sg)),
\label{rhom}
\end{eqnarray}
where
\begin{equation}\label{alpha-p}
  \a_p(X)=\max\bigl(\,d_{2p},d_{2p+2}\,\bigr),\quad d_k=\norm{X^k}^{1/k}
\end{equation}
and $p$ is any positive integer satisfying the constraint $m+1\ge p(p-1)$.
In addition, it is straightforward to verify that
\begin{equation}
 \norm{\sinc(s^{-1}A^\sg)-\rt_m(s^{-1}A^\sg)} \le \sum_{j=m+1}^{\infty}\dfrac{\a_p(s^{-1}A^\sg)^{2j}}{(2j+1)!}
                    \le
\rho_m(\a_p(s^{-1}A^\sg)).
\end{equation}
Similarly the forward errors of the approximations of
$\cosh$, and $\sinch$ by Taylor polynomials have exactly
the same bound $\rho_m$.

Next we analyze the computational cost and determine how to
choose the scaling parameter and the degree of Taylor polynomial.
Define
\begin{equation}\label{thm}
 \th_m=\max\{\,\th : \rho_m(\th)\le\tol\,\}.
\end{equation}
Thus given $m$ and $p$ if $s$ is chosen so that $s^{-1}\a_p(A^\sg)\le\th_m$,
then the inequality $\rho_m(\a_p(s^{-1}A^\sg))\le\tol$ will be satisfied and
therefore the absolute forward error will be bounded by $\tol$.
Table~\ref{table.thetavals} lists selected values of $\th_m$ for
$\tol=2^{-10}$ (half precision), $\tol=2^{-24}$ (single precision),
and $\tol=2^{-53}$ (double precision). These values were determined as
described in \cite[App.]{hial10}.
For each $m$, the optimal value of the scaling parameter $s$ is given
by $s=\max(\lceil\a_p(A^\sg)/\th_m\rceil,1)$. The computational cost of
evaluating $T_s$ in view of Code Fragment \ref{code.frag1} is $2\sg ms$
matrix--matrix multiplications of the form $AB$. That is,
$2\sg n_0ms$ \matvec\ products since $B$ has $n_0$ columns.
By \eqref{ev.od.sum}, $U_{s-1}$ is obtained with a negligible cost.
$\sinc(A)B$ can be then recovered by a single invocation of
Code Fragment \ref{code.frag1} for $V=U_{s-1}$ and $\gamma=2k+1$; this
requires only $2\sg n_0m$ \matvec\ products. After that one
multiplication is needed to recover $\sin(A)B$ from $\sinc(A)B$;
that is $n_0$ \matvec\ products.
We build our cost analysis on an assumption that the output of our
\alg\ is $\cos(A^\sg)B$ and $\sinc(A^\sg)B$. Note that when $\sg=1/2$,
$\sin(A^\sg)B$ cannot be obtained without computing $\sqrtA$.
Thus the total cost is
\begin{equation}\label{Cm}
 2\sg n_0m(s+1) = 2\sg n_0m\bigl(\max(\lceil\a_p(A^\sg)/\th_m\rceil,1)+1\bigr)
\end{equation}
\matvec\ products.
We observe that this quantity tends to be decreasing as $m$ increases
though the decreasing is not necessarily monotonic.
The sequence $\{m/\th_m\}$ is strictly decreasing while the sequence
$\{\a_p(X)\}$ has a generally nonincreasing trend for any $X$.
Thus the larger is $m$, the less the cost.
However, a large value of $m$ could lead to unstable calculation of
Taylor polynomials $r_m(A^\sg)B$ for large $\norm{A^\sg}$ in \fpa.
Thus we impose a limit $m_{\max}$ on $m$ and seek $m_*$ that minimizes
the computational cost over all $p$ such that $p(p-1)\le m_{\max}+1$.
For the moment we drop the $\max$ in
\eqref{Cm}, whose purpose is simply to cater for nilpotent $A^\sg$ with
$A^{\sg j} = 0$ for $j \ge 2p$. Moreover, we remove constant terms since
they essentially don't effect the optimization for the value of $m_*$.
Thus we consider the sequence
$$
    C_m(A^\sg) = m\lceil\a_p(A^\sg) /\th_m \rceil
$$
to be minimized subject to some constraints.
Note that $\norm{A^\sg}\ge d_2 \ge d_{2k}$ in \eqref{alpha-p} for all $k\ge 1$ and so
\begin{equation}\label{ineq}
 \norm{A^\sg}\ge\a_1(A^\sg)=d_2=\norm{A^{2\sg}}^{1/2} \ge \a_p(A^\sg)
\end{equation}
for all $p\ge1$.
Hence we don't need to consider the case $p = 1$ when minimizing
$C_m(A^\sg)$ since $C_m(A^\sg) \le m\lceil\a_1(A^\sg) /\th_m \rceil$ . Let
$\pmax$ denote the largest positive integer $p$ such that $p(p-1)\le
m_{\max}+1$.
Let $m_*$ be the smallest value of $m$ at which the minimum
\begin{eqnarray}
   \qquad
    \Cmstar(A^\sg) &=& \min\bigl\{\,m\lceil\a_p(A^\sg) /\th_m \rceil :
    \mbox{$2 \le p \le \pmax$,~~$p(p-1)-1\le m \le \mmax$}
\,\bigr\},
   \label{Cmopt}
\end{eqnarray}
is attained \cite[Eq.~(3.11)]{alhi11}.
The optimal scaling parameter then is
$$
s = \max(\Cmstar(A^\sg)/m_*,1).
$$
Our experience and observation indicate that $\pmax = 5$ and $\mmax = 25$ are
appropriate choices for our \alg. However the \alg\ supports user-specified
values of $\pmax$ and $\mmax$.

The forward error analysis and cost analysis
are valid for any matrix norm, but it is most convenient to use the 1-norm
since it is easy to be efficiently estimated using the block 1-norm
estimation \alg\ of Higham and Tisseur \cite{hiti00n}.
We estimate the quantities $d_k = \normi{A^{\sg k}}^{1/k}$,
where $k$ is even as defined in \eqref{alpha-p},
which are required to form $\a_p(A^\sg)$.
The \alg\ of Higham and Tisseur estimates $\normi{A^{\sg k}}$ via about
two actions of $A^{\sg k}$ and two actions of $(A^*)^{\sg k}$, all on
matrices of $\ell$ columns, where the positive integer $\ell$ is a
parameter (typically set to 1 or 2).
The number $\sg k$ is a positive integer since $k$ is even, so fractional
powers of $A$ is completely avoided.
Therefore obtaining $\a_p(A^\sg)$ for $p=2\colon p_{\max}$
costs approximately
\begin{equation}
    8\sg\ell\sum_{p=2}^{p_{\max}+1}p = 4\sg@\ell@ p_{\max}(p_{\max}+3)
  \label{capprox}
\end{equation}
\matvec\ products. Thus in view of \eqref{Cm} if it happens that
$$
2\sg n_0\mmax\bigl(\normi{A}^\sg/\th_{\mmax}+1\bigr)\le 4\sg@\ell@ p_{\max}(p_{\max}+3),
$$
or equivalently
\begin{equation}
   \normi{A}^\sg \le \th_{\mmax}\left(\frac{2\ell}{n_0\mmax}
                            p_{\max}(p_{\max}+3)-1\right)
   \label{usenormA}
\end{equation}
then the computational cost of evaluating $T_s$ and
$\rt_m(s^{-1}A^\sg)U_{s-1}$ with $m$
determined by using $\normi{A}^\sg$ or $\normi{A^{2\sg}}^{1/2}$
in place of $\a_p(A^\sg)$ in \eqref{Cmopt} is no larger than the cost
\eqref{capprox} of computing the sequence $\{\a_p(A^\sg)\}$.
Thus we should certainly use $\normi{A}^\sg$ if $\sg=1$ or $\normi{A}^{1/2}$
if $\sg=1/2$ in place of $\a_p(A^\sg)$ for each $p$
in light of the inequalities in \eqref{ineq}.

In the case $\sg=1$, we still have another chance to avoid estimating
$\a_p(A)$ for $p>2$.
If the inequality \eqref{usenormA} is unsatisfied,
the middle bound $d_2$ in \eqref{ineq} can be estimated and its actual cost,
$\nu$ \matvec\ products, can be counted. We check again if the bound
\begin{equation}
   d_2 \le \th_{\mmax}\left(\frac{2\ell}{n_0\mmax}
                            p_{\max}(p_{\max}+3)-\nu-1\right)
   \label{usenormd2}
\end{equation}
holds. We sum up our analysis for determining the parameters $m_*$ and $s$
in the following code.
\begin{fragment}[\mbox{$[m_*,s,\Theta_\sg] = \parameters(A,\sg,\tol)$}]\label{code.frag0}
This code determines $m_*$ and $s$
given $A$, $\sg$, $\tol$, $\mmax$, and $\pmax$.
Let $\Theta_\sg$ denote the number of the actual \matvec\ products
needed to estimate the sequence $\{\a_p(A^\sg)\}$.
\end{fragment}

\begin{code}
if \= \eqref{usenormA} is satisfied\\
\>$m_* = \argmin_{1\le m\le \mmax} m \lceil \normi{A}^\sg/\th_m \rceil$\\
\>$s = \lceil \normi{A}^\sg/\th_{m_*} \rceil$\\
\>goto line \ref{end}\\
end\\
if $\sg=1$\\
\> Compute $d_2$\\
\> if \= \eqref{usenormd2} is satisfied\\
\>\> $m_* = \argmin_{1\le m\le \mmax} m \lceil d_2/\th_m \rceil$\\
\>\> $s = \lceil d_2/\th_{m_*} \rceil$\\
\>\> goto line \ref{end}\\
\>end\\
end\\
Let $m_*$ be the smallest $m$ achieving the minimum in \eqref{Cmopt}.\\
$s = \max\bigl( \Cmstar(A^\sg)/m_*, 1\bigr)$\\
\label{end}
end
\end{code}
As explained in \cite[sect.~3]{alhi11},
if we wish to compute $f(tA^\sg)B$ for several values of $t$, we need not
invoke Code Fragment \ref{code.frag0} for each $t^{1/\sg}A$.
The trick is that since $\a_p(tA^\sg) = |t|@ \a_p(A^\sg)$, we can precompute the matrix
$S\in\R^{(\pmax-1)\times \mmax}$ given by
\begin{equation}
   S_{pm} = \mycases{ \dfrac{\a_{p}(A^\sg)}{\th_m}, & $2 \le p \le \pmax$,~
                                                $p(p-1)-1\le m \le \mmax$,\cr
                          0, & otherwise\cr}
   \label{spm}
\end{equation}
and then for each $t$ obtain $\Cmstar(tA^\sg)$ as the smallest nonzero
element in the matrix $\lceil |t|S \rceil \diag(1,2,\dots,\mmax)$,
where $m_*$ is the column index of the smallest element.
The benefit of basing the selection of the scaling parameter
on $\a_p(A)$ instead of $\norm{A}$ is that $\a_p(A)$ can be much smaller
than $\norm{A}$ for highly nonnormal matrices.

\begin{table}
%
\caption{Selected constants $\th_m$ for $\tol=2^{-10}$ \parens{half}, $\tol = 2^{-24}$
        \parens{single}, and $\tol = 2^{-53}$
        \parens{double}.
        }
\label{table.thetavals}
\setlength{\tabcolsep}{4.5pt}
\begin{center}\footnotesize
\begin{tabular}{c|*{12}{c}}
$2m$ &  6 & 10 & 14 & 18 & 22 & 26 & 30 & 34 & 38 & 42 & 46 & 50\\\hline\mystrut{8pt}
 half &   1.6e0 &   3.0e0 &   4.4e0 &   5.8e0 &   7.3e0 &   8.8e0 &   1.0e1 &   1.2e1 &   1.3e1 &   1.5e1 &   1.6e1 &   1.8e1\\
single &   1.8e0 &   4.2e0 &   6.9e0 &   9.7e0 &   1.3e1 &   1.5e1 &   1.8e1 &   2.1e1 &   2.4e1 &   2.7e1 &   3.0e1 &   3.3e1\\
double &   3.8e-2 &   2.5e-1 &   6.8e-1 &   1.3e0 &   2.1e0 &   3.0e0 &   4.1e0 &   5.1e0 &   6.3e0 &   7.5e0 &   8.7e0 &   1.0e1\\

\end{tabular}
\end{center}

\end{table}
\section{Preprocessing and termination criterion}\label{sec.prep}
In this section we discuss several strategies to improve the \alg\
stability and reduce its computational cost. The algorithmic scaling parameter $s$
plays an important role that the smaller the $s$ the better the
stability of the \alg\, in general, and the lower the computational cost.
That why we rely on $\a_p(A)$ instead of merely using $\norm{A}$ to produce
the scaling parameter.
Al-Mohy and Higham \cite[sect.~3.1]{alhi11} proposed an argument reduction
and a termination criterion. They have found empirically that the shift
$\mu=n^{-1}\trace(A)$ \cite[Thm.~4.21]{high:FM}
that minimizes the Frobenius norm of the matrix $\At=A-\mu I$ leads to
smaller values of $\a_p(\At)$ than $\a_p(A)$. We use this shift here
if the required outputs are $\cos(A)B$ and $\sin(A)B$ or
$\cosh(A)B$ and $\sinh(A)B$.
There are cases where shifting is impossible to recover.
This happens when the required outputs include $\sinc(A)B$, $\sinch(A)B$,
or any form of $f(\sqrtA)B$.

We can recover the original cosine and sine of $A$ from the computed
cosine and sine of $\At$ using the formulas
\begin{equation}
 \label{shift.cos.sin}
  \cos A = \cos\mu \cos \At - \sin\mu \sin \At,\quad
  \sin A = \cos\mu \sin \At + \sin\mu \cos \At.
\end{equation}
The functions
$\cosh A$ and $\sinh A$ have analogous formulas containing $\cosh\mu$
and $\sinh\mu$, which could \emph{overflow} for
large enough $|\Re(\mu)|$. Same problem arises for $\cos\mu$ and
$\sin\mu$ if $|\Im(\mu)|$ is large enough.
Al-Mohy and Higham successfully overcome this problem in their \alg\
for the matrix
exponential by undoing the effect of
the scaled shift right after the inner loop of \cite[Alg.~3.2]{alhi11}.
It is possible to do so for trigonometric and hyperbolic matrix functions.
We can undo the effect of the scaled shift in $\cos(s^{-1}\At)T_{k-1}$ for
each $k$ in the recurrence \eqref{rec.cos} using the formula in
\eqref{shift.cos.sin}, which requires $\sin(s^{-1}\At)T_{k-1}$.
The next code shows how $\sin(s^{-1}\At)T_{k-1}$ can be formed
using the already generated power actions, $\At^{2k}B$.

\begin{fragment}\label{code.frag3}
Given $\At=A-\mu I\in\C^{\nbyn}$, $B\in\C^{n\times n_0}$, and a suitable
chosen scaling parameter $s$, this code returns
$C=r_m(s^{-1}A)B\approx\cos(s^{-1}A)B$.
\end{fragment}
\begin{code}
$V = B$, $Z = B$\\
for \= $k=1:m$\\
\label{line.gamma.q}
\> $\beta = 2k$, $\gamma = 2k-1$, $q=1/(2k+1)$ \\
\> $B = \At B$\\
\> $B = (\At B)\left(s^2\beta\gamma\right)^{-1}$\\
\> $V = V + (-1)^kB$\\
\> $Z = Z + (-1)^kqB$\\
end\\
$C = \cos(\mu/s)V-s^{-1}\sin(\mu/s)\At Z$
\end{code}

The recovery of $\sin(s^{-1}A)U_{s-1}$ (recall \eqref{sin.U})
can be obtained by a single
execution of Code Fragment \ref{code.frag3} for $V=Z=U_{s-1}$ after
setting $\gamma=2k+1$ and $q=2k+1$ in line \ref{line.gamma.q}. Thus
$\sin(s^{-1}A)U_{s-1}\approx s^{-1}\cos(\mu/s)\At V + \sin(\mu/s)Z$.
Comparing Code Fragment \ref{code.frag1} with Code Fragment
\ref{code.frag3} assuming the same scaling parameter $s$,
undoing the shift requires $n_0$ \matvec\ products for each
$k=1\colon s$ bringing the total of the extra cost to
$(s+1)n_0$ \matvec\ products: $sn_0$ for $T_s$ and $n_0$ to recover
$\sin(A)B$ from $\sin(s^{-1}\At)U_{s-1}$ using \eqref{sin.U} and the
formula in \eqref{shift.cos.sin}. However, the scaling parameter
$s$ selected based on $\At$ is potentially smaller than that
selected based on $A$ making the overall cost of the \alg\ potentially
smaller.

For the early termination of the evaluation of Taylor polynomials,
we use the criterion proposed by Al-Mohy and Higham
\cite[Eq.~(3.15)]{alhi11} implemented in line \ref{line.a12} of \Alg\
\ref{alg.funmv} below.
\section{Algorithm}\label{sec3}

In this section we write in details our \alg\ for computing the
trigonometric and hyperbolic matrix functions of the forms:
$f(tA)B$ and $f(t\sqrtA)B$.

\begin{algorithm}[\t{$[C,S]$ = funmv($t,A,B,\tol$,option)}]\label{alg.funmv}
Given $t\in\C$, $A\in\C^{\nbyn}$, $B\in\C^{n\times n_0}$, and
a tolerance $\tol$, this \alg\
computes $C$ and $S$ for any chosen option of the table.
The parameters $\sg$, $k_0$, and \shift are set to their
corresponding values of the last column depending on the chosen case.
\end{algorithm}
\vspace{1.0mm}
\begin{center}
\begin{tabular}{|c |m{3cm} m{3cm} |c| }
 \hline
option & \qquad\qquad outputs    &                             & $(\sg,k_0,\shift)$\\[4pt]\hline\mystrut{10pt}
$1$& $C\approx\cos(tA)B$       &  $S\approx\sin(tA)B$        & $(1,1,1)$         \\[4pt]\hline\mystrut{10pt}
$2$& $C\approx\cosh(tA)B$      &  $S\approx\sinh(tA)B$       & $(1,0,1)$         \\[4pt]\hline\mystrut{10pt}
$3$& $C\approx\cos(tA)B$       &  $S\approx\sinc(tA)B$       & $(1,1,0)$         \\[4pt]\hline\mystrut{10pt}
$4$& $C\approx\cosh(tA)B$      &  $S\approx\sinch(tA)B$      & $(1,0,0)$         \\[4pt]\hline\mystrut{10pt}
$5$& $C\approx\cos(t\sqrtA)B$  & $S\approx\sinc(t\sqrtA)B$   & $(\half,1,0)$     \\[4pt]\hline\mystrut{10pt}
$6$& $C\approx\cosh(t\sqrtA)B$ & $S\approx\sinch(t\sqrtA)B$  & $(\half,0,0)$     \\[4pt]\hline
\end{tabular}
\\[10pt]
\end{center}

\begin{code}
if \shift, $\mu = \trace(A)/n$, $A = A-\mu I$, end\\
if \= $t\normi{A} = 0$\\
   \> $m_* = 0$, $s = 1$ ~~\% The case $tA=0$.\\
else\\
\label{line.code.frag}
 \> $[m_*,s,\Theta_\sg] = \parameters(t^{1/\sg}A,\sg,\tol)$  \% Code Fragment~\ref{code.frag0}\\
end\\
$\undoin = 0$, $\undout = 0$ ~~\% undo shifting inside or outside the loop.\\
if \= option 1 and $|\Im(t\mu)|>0$\\
\> $\phi_1=\cos(t\mu/s)$, $\phi_2=\sin(t\mu/s)$, $\undoin = 1$\\
elseif option 1 and $t\mu\in\R\backslash\{0\}$\\
\> $\phi_1=\cos(t\mu)$, $\phi_2=\sin(t\mu)$, $\undout = 1$\\
elseif option 2 and $|\Re(t\mu)|>0$\\
\> $\phi_1=\cosh(t\mu/s)$, $\phi_2=\sinh(t\mu/s)$, $\undoin = 1$\\
elseif option 2 and $t\mu\in\C\backslash\R$\\
\> $\phi_1=\cosh(t\mu)$, $\phi_2=\sinh(t\mu)$, $\undout = 1$\\
end\\
$T_0 = 0$\\
if $2|s$, $T_0=B/2$, end\\
$U=T_0$, $T_1=B$\\

for \= $i=1:s+1$\\
\> if \= $i=s+1$\\
\>\> $U=2U$, $T_1=U$\\
\> end\\
\> $V=T_1$, $Z=T_1$, $B=T_1$\\
\> $c_1 = \normo{B}$\\
\>  for \= $k = 1:m_*$\\
\> \>    $\beta=2k$\\
\> \>    if $i\le s$, $\gamma=\beta-1$, $q=1/(\beta+1)$
        else $\gamma=\beta+1$, $q=\gamma$, end\\
\label{line.sqrt}
\> \>   if $\sg=1$, $B = AB$, end\\
\> \>   $B = (AB)\bigl((t/s)^2/(\beta\gamma)\bigr)$\\
\> \>   $c_2 = \normo{B}$\\
\label{line.V}
\> \>   $V = V + (-1)^{k_0k}B$\\
\>\>    if \undoin, $Z = Z + ((-1)^{k_0k}q)B$, end\\

\label{line.a12}
\> \>   if \= $c_1 + c_2 \le \tol@\normo{V}$, 
        break, 
        end\\
\> \>   $c_1 = c_2$\\
\>  end\\
\> if \= \undoin\\
\>\>   if \= $i \le s$\\
\>\>\>    $V = V\phi_1 + A(Z((-1)^{k_0}t\phi_2/s))$\\
\>\>   else\\
 \>\>\>   $V = A(V(t\phi_1/s)) + Z\phi_2$\\
\>\>   end\\
\> end\\
\> if \=$i = 1$, $T_2 = V$, elseif $i \le s$, $T_2 = 2V - T_0$, end
 ~~ \% using \eqref{rec.cos}.\\
\>    if \=$i \le s-1$ and $(2|s \xor 2|i)$\\
\>\>     $U = U + T_2$ ~~ \% using \eqref{ev.od.sum}.\\
\>    end\\
\>    $T_0 = T_1$,  $T_1 = T_2$\\
end\\
$C = T_2$\\
if \= \undoin\\
\> $S = V$\\
elseif option 1 or option 2\\
\> $S = A(V(t/s))$\\
else\\
\> $S = V/s$\\
end\\
if \= \undout\\
\>   $C = \phi_1C + ((-1)^{k_0}\phi_2)S$\\
\>   $S = \phi_1S + \phi_2T_2$\\
end
\end{code}

Due to the stopping criterion in line \ref{line.a12}, assume that
the inner loop is terminated when $k$ takes values $m_i$, $i=1\colon s+1$.
Thus the total cost of the \alg\ is
$$
2\sg n_0\sum_{i=1}^{s+1}m_i +(\undoin)n_0(s+1)+(\shift-\undoin)n_0+\Theta_\sg
$$
\matvec\ multiplications. Since $m_i$ and $\Theta_\sg$ are bounded by $m_*$  and
\eqref{capprox}, \resp, an upper bound of the computational cost of the \alg\
can be obtained after the execution of Code Fragment \ref{code.frag0}
in line \ref{line.code.frag}.
This advantage allows users to estimate the overhead of the \alg.
When $\sg=1/2$, the \alg\ saves about 50 percent of the computational
cost comparing to the other options. Therefore it is better not to provide $\sqrtA$
even if it is easy to evaluate. As an example, take
$A=\diag(1,2,\cdots,100)$, $b=[1,1,\cdots,1]^T$, and $t=1$.
Executing \t{funmv($t,A,B$)} (option 5) requires 51 \matvec\ products whereas
\t{funmv($t,A^{1/2},B$)} (option 3) requires 102.
Note that it is possible to obtain $\sin(tA)B$ and $\sinh(tA)B$ in options 3 and 4,
\resp. However this is impossible in options 5 and 6 because of the absence
of $\sqrtA$. The present of shifting in options 1 and 2 makes it impossible
to obtain $\sinc(tA)B$ and $\sinch(tA)B$ as we pointed out in the previous
section.


\section{Numerical experiments}\label{sec4}
In this section we give some numerical tests to illustrate the accuracy
and efficiency of \Alg~\ref{alg.funmv}. We use MATLAB$^\circledR$ R2015a on a machine with
Core i7. The experiments involve the following \alg s:
\begin{figure}
  \centering
  \Myfig{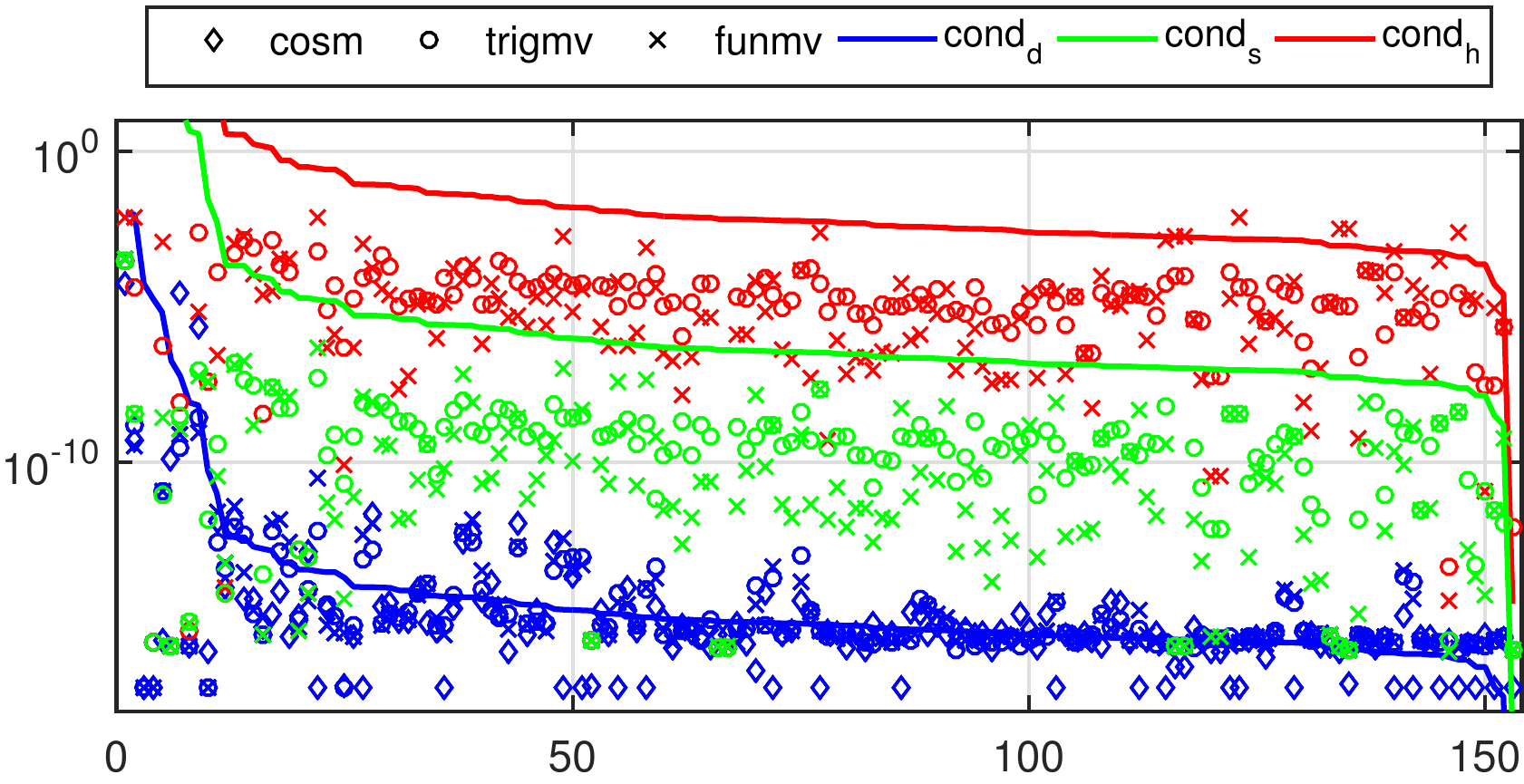}
  \caption{Experiment \ref{ex1}. Normwise relative errors
  in computing $\cos(A)b$ using different precisions.
  $\cond_d$, $\cond_s$, and $\cond_h$ represent
  $\cond(\cos,A)$ multiplied by
  $2^{-53}$, $2^{-24}$, and $2^{-10}$, \resp.
          }
  \label{fig.test1}
\end{figure}
\begin{figure}
  \centering
  \Myfig{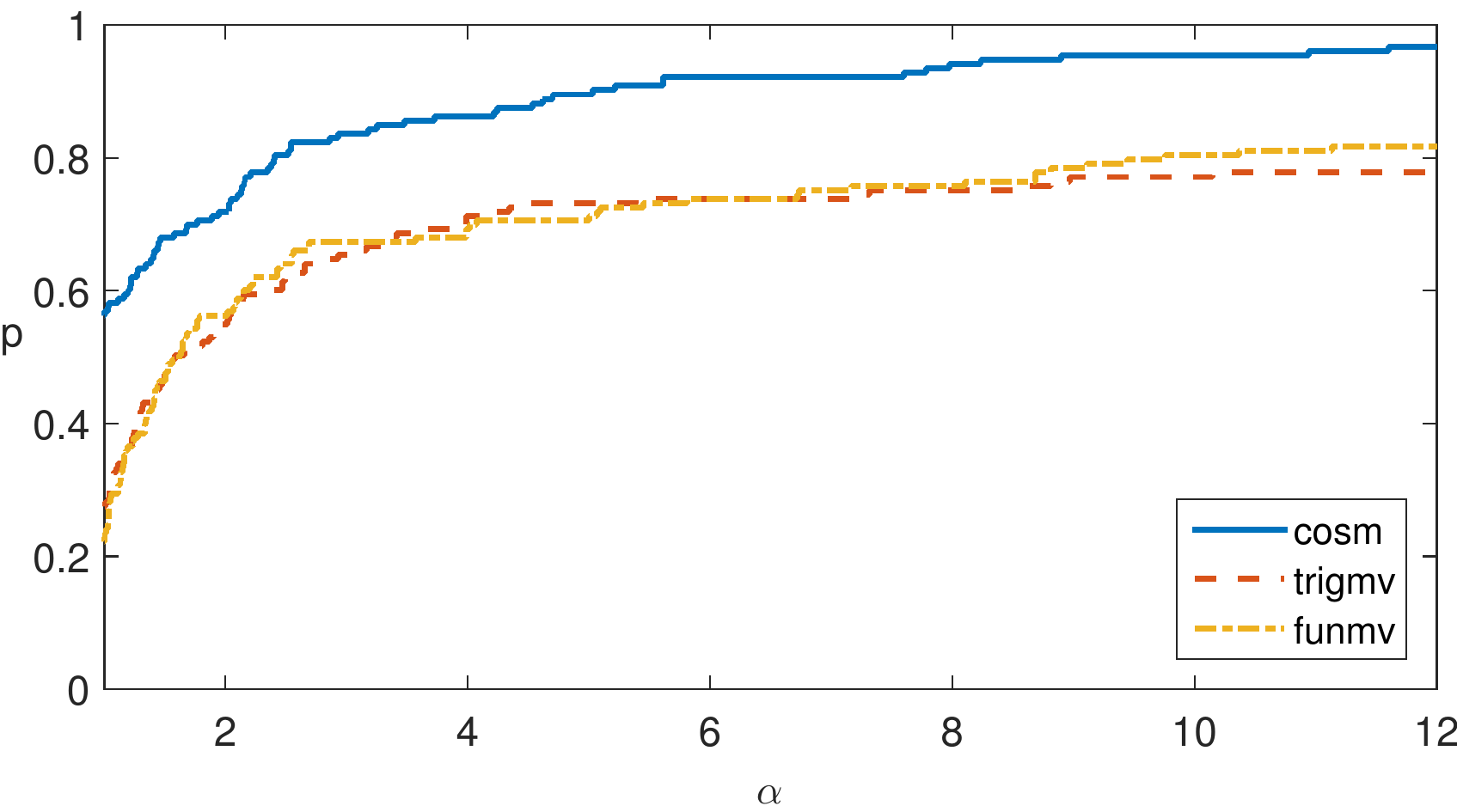}
  \caption{Double precision data of Figure \ref{fig.test1}
           presented as a performance profile.
          }
  \label{fig.pef}
\end{figure}
\begin{enumerate}
  \item \t{funmv}: the MATLAB code of \Alg~\ref{alg.funmv},
  \item \t{trigmv} and \t{trighmv}: MATLAB codes implementing
   the recently authored \alg\ by Higham and Kandolf
  \cite[Alg.~3.2]{hika17}. \t{trigmv} returns the actions
  $\cos(A)b$ and $\sin(A)b$ while \t{trighmv} returns the actions
  $\cosh(A)b$ and $\sinh(A)b$. The codes are available in
   \url{https://bitbucket.org/kandolfp/trigmv}
   \item \t{expmv}: MATLAB code for the \alg\ of Al-Mohy and Higham
   \cite[Alg.~3.2]{alhi11} that compute the action of the matrix
   exponential $e^AB$. The code is available in
   \url{https://github.com/higham/expmv}.
  \item \t{cosm} and \t{sinm}: \cite[Alg's~4.2 \& 5.2]{ahr15}
  of Al-Mohy, Higham, and Relton
  for explicitly computing $\cos A$ and $\sin A$, \resp.
  The multiplication by $b$ follows to obtain $\cos(A)b$
  or $\sin(A)b$. The MATLAB codes of the \alg s
  are available in \url{https://github.com/sdrelton/cosm_sinm}.
\end{enumerate}
\begin{figure}[t]
  \centering
  \Myfig{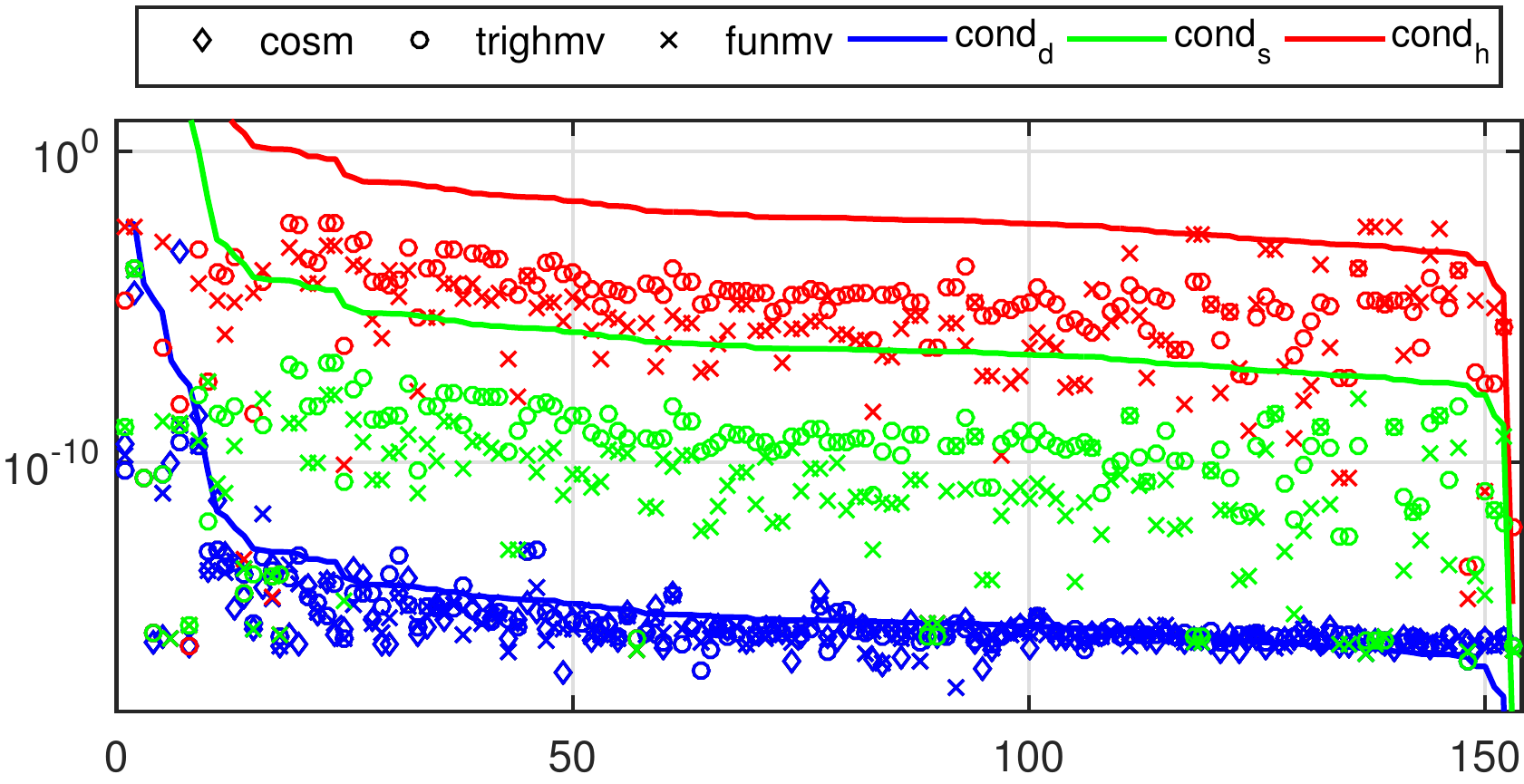}
  \caption{Experiment \ref{ex1}. Normwise relative errors
  in computing $\cosh(A)b$ using different precisions.
  $\cond_d$, $\cond_s$, and $\cond_h$ represent
  $\cond(\cosh,A)$ multiplied by
  $2^{-53}$, $2^{-24}$, and $2^{-10}$, \resp.
          }
  \label{fig.test2}
\end{figure}
\begin{figure}[t]
  \centering
  \Myfig{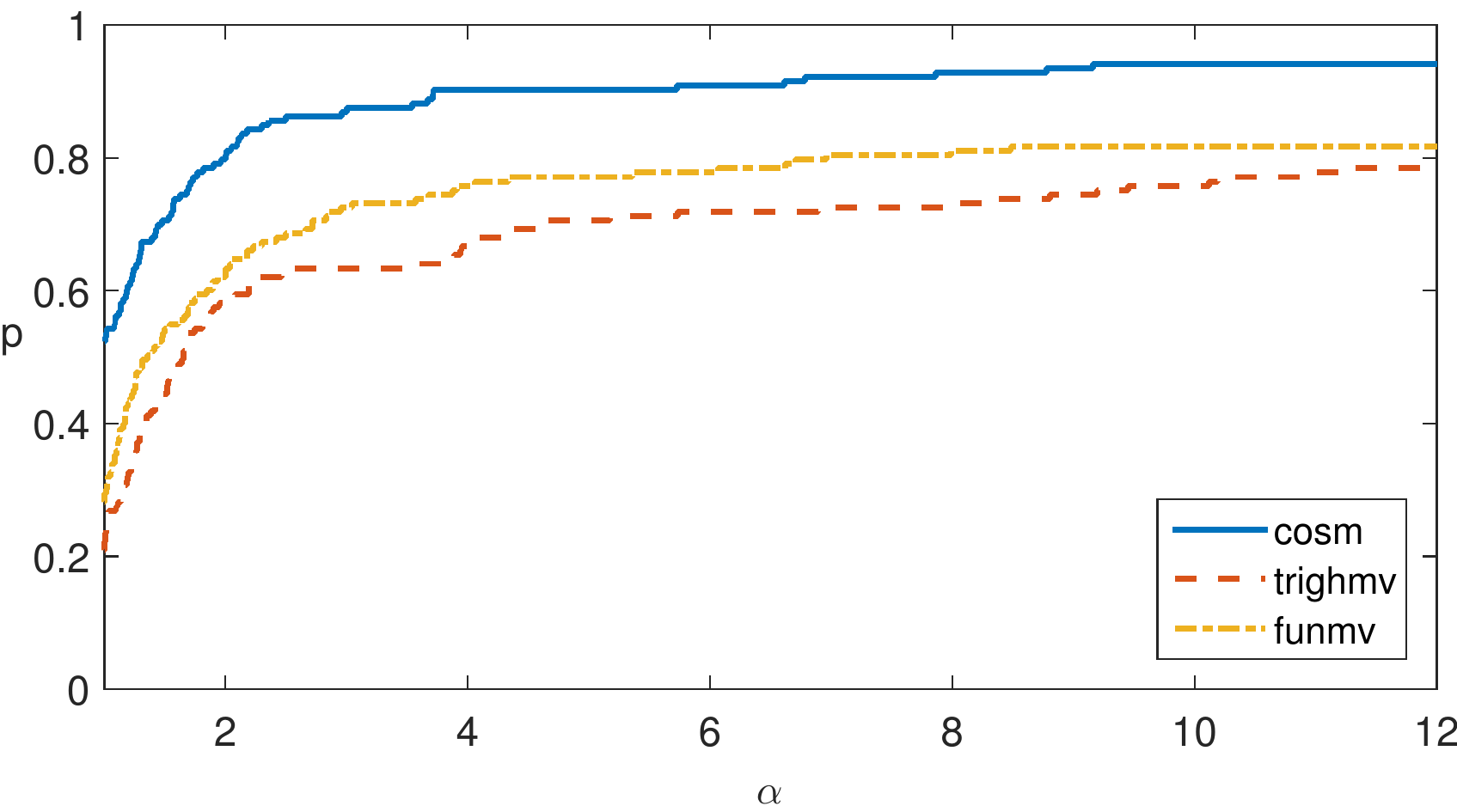}
  \caption{Double precision data of Figure \ref{fig.test2}
           presented as a performance profile.
          }
  \label{fig.pef2}
\end{figure}
\begin{example}
\item\label{ex1}
In this experiment we test the stability of \t{funmv} (option 1) comparing with
\t{trigmv} and \t{cosm}. We use the test matrices described in
\cite[sect.~6]{alhi09a} and used also in \cite[sect.~6]{alhi11}.
For each matrix $A$ of these test matrices, a vector $b$ is
randomly generated. We approximate $x:=\cos(A)b$ by $\xhat$ using \t{funmv} and
\t{trigmv} with the tolerances of half, single, and double precisions. The approximation
of $x$ by \t{cosm} is carried out in double precision since the \alg\
is only intended for that.
The ``exact" $x$ is computed
at 100 digit precision with the Symbolic Math Toolbox. The relative forward
errors $\normt{x-\xhat}/\normt{x}$ for each tolerance is plotted in
Figure \ref{fig.test1}, where the solid lines represent the condition
number of the matrix cosine $\cond(\cos,A)$ multiplied by the associate
tolerance $\tol$ sorted in a descending order.
The condition number with respect to Frobenius norm is estimated
using the code \t{funm\_condest\_fro}
from the Matrix Function Toolbox \cite{high-mft}.

Figure \ref{fig.pef} displays a performance profile
for the double precision data plotted in Figure \ref{fig.test1}
which includes the data of \t{cosm}.
For each method, the parameter $p$ is the proportion of problems
in which the error is within a factor of $\a$ of the smallest
error over all methods.
The experiment reveals that our \alg\ behaves as stable as the existing
\alg s. The performance profile shows that \t{cosm} outperforms the
other methods while \t{funmv} and \t{trigmv} have similar behavior.

We repeat the experiment for $\cosh(A)b$ using \t{funmv} (option 2),
\t{trighmv}, and \t{cosm} with argument $\i A$.
The results are reported in Figure \ref{fig.test2} and Figure
\ref{fig.pef2}. Both methods behave in a stable manner but \t{funmv}
outperforms \t{trighmv} in view of the performance profile.

The figures corresponding to $\sin(A)b$ and $\sinh(A)b$ are similar
to those of $\cos(A)b$ and $\cosh(A)b$, \resp; that why we don't report
them here.
\item \label{ex2}
In this experiment we compute $\cos(A)b$ for large and sparse matrices.
We compare \t{funmv} (option 1) with \t{trigmv} in terms of CPU time,
\matvec\ products, and relative forward errors in 1-norm.
We use \t{cosm} to compute the reference solution in double precision.
The test matrices are prescribed in \cite[Example 4.2]{hika17} and
\cite[Experiment 5]{alhi11}. The first three matrices of Table
\ref{table.ex2} belong to the Harwell–-Boeing collection and are
obtained from the University of Florida Sparse Matrix Collection \cite{dahu11}.
The matrix \t{triw} and \t{poisson} are from the MATLAB
\t{gallery}. The matrices and problem details are
\begin{itemize}
  \item \t{orani678} (nonsymmetric), $n=2529$, $b=[1,1,\cdots,1]^T$;
  \item \t{bcspwr10} (symmetric), $n=5300$, $b=[1,0,\cdots,0,1]^T$;
  \item \verb"gr_30_30", $n=900$, $b=[1,1,\cdots,1]^T$;
  \item \t{triw} denotes \t{-gallery('triw',2000,4)} (upper triangular
  with $-1$ in the main diagonal and $-4$ elsewhere), $n=2000$,
  $b=[\cos 1, \cos 2,\cdots,\cos n]^T$;
  \item \t{poisson} denotes \t{-gallery('poisson',99)}
  (symmetric negative definite),
  $n=9801$,
  $b=[\cos 1, \cos 2,\cdots,\cos n]^T$. This matrix arises
  from a finite difference discretization of the two--dimensional
  Laplacian in the unit square.
\end{itemize}

\begin{table}[t]  


\caption{\Experiment\ \ref{ex2}:
         \speed\ denotes time for method divided by
         time for \t{funmv}.}

\label{table.ex2}

\begin{center} \footnotesize

(a) Double precision\\
\begin{tabular}{cc|*{3}{c}|*{3}{c}|*{1}{c}}
\\\hline
&& \multicolumn{3}{c|}{\t{funmv} } &
   \multicolumn{3}{c|}{\t{trigmv}} &
   \multicolumn{1}{c}{\t{cosm}}\\
& $t$ & \speed  & \mv  & Error  & \speed & \mv   & Error & \speed\\
\hline
\mystrut{8pt}
\t{orani678}    &  100 &   1 &     1111 &   6.0e-15  &   1.4 &     2024 &   4.5e-15&   9.9e2 \\
\t{bcspwr10}    &  10  &   1 &      379 &   3.8e-14  &   1.7 &      618 &   3.8e-14&   2.5e3 \\
\verb"gr_30_30" &  2   &   1 &      133 &   6.1e-14  &   1.3 &      188 &   7.8e-14&   3.2e2 \\
\t{triw}        &  10  &   1 &    27005 &   7.1e-14  &   1.2 &    56560 &   1.4e-13&   2.2e-1\\
\t{poisson}     &  500 &   1 &     9757 &   4.0e-13  &   2.2 &    19036 &   2.2e-13&   1.0e3
\end{tabular}
\vspace{0.3cm}
\\
(b) Single precision\\
\begin{tabular}{cc|*{3}{c}|*{3}{c}|*{1}{c}}
\\\hline
&& \multicolumn{3}{c|}{\t{funmv} } &
   \multicolumn{3}{c|}{\t{trigmv}} &
   \multicolumn{1}{c}{\t{cosm}}\\
& $t$ & \speed  & \mv  & Error  & \speed & \mv   & Error & \speed\\
\hline
\mystrut{8pt}
\t{orani678}    &  100 &   1 &   719 &   2.1e-9  &   1.5 &  1224  &   4.1e-8 &   1.8e3 \\
\t{bcspwr10}    &  10  &   1 &   265 &   3.0e-10 &   1.6 &  402   &   4.7e-10&   3.7e3 \\
\verb"gr_30_30" &  2   &   1 &    97 &   3.8e-9  &   1.2 &  136   &   5.5e-9 &   5.1e2 \\
\t{triw}        &  10  &   1 & 13011 &   8.2e-13 &   1.1 &  26708 &   5.1e-9 &   4.4e-1\\
\t{poisson}     &  500 &   1 & 6415  &   1.3e-8  &   2.2 &    12436 &   2.5e-7&   1.5e3
\end{tabular}
\vspace{0.3cm}
\\
(c) Half precision\\
\begin{tabular}{cc|*{3}{c}|*{3}{c}|*{1}{c}}
\\\hline
&& \multicolumn{3}{c|}{\t{funmv} } &
   \multicolumn{3}{c|}{\t{trigmv}} &
   \multicolumn{1}{c}{\t{cosm}}\\
& $t$ & \speed  & \mv  & Error  & \speed & \mv   & Error & \speed\\
\hline
\mystrut{8pt}
\t{orani678}    &  100 &   1 &      551 &   1.6e-4  &   1.3 &    848 &   1.9e-3&   2.3e3 \\
\t{bcspwr10}    &  10  &   1 &      215 &   7.3e-6  &   1.5 &    324 &   1.3e-5&   4.4e3 \\
\verb"gr_30_30" &  2   &   1 &       93 &   3.1e-4  &   1.3 &    108 &   5.4e-6&   4.9e2 \\
\t{triw}        &  10  &   1 &     7381 &   1.1e-6  &   1.1 &  15028 &   5.4e-5&   7.6e-1\\
\t{poisson}     &  500 &   1 &     5223 &   2.7e-4  &   2.1 &   9810 &   4.2e-4&   1.9e3
\end{tabular}
\end{center}

\end{table}
The results are shown in Table \ref{table.ex2}. The three blocks of the
table display the computations with different tolerances that represent
double, single, and half precisions. The symbol \speed\ denotes
CPU time for method divided by CPU time for \t{funmv} and the symbol
\mv\ denotes the number of matrix--vector products required by each methods.
Obviously \t{funmv} proves superiority. It does  outperform \t{trigmv} in terms
of CPU running time and computational cost. the number of matrix--vector
products of \t{funmv} is about the half of that of \t{trigmv} for most cases.
No wonder since \t{trigmv} requires the action of the matrix exponential
on a matrix of two columns---namely
$B=[b,b]/2$---to yield $\cos(tA)b$.

\item \label{ex3}
In this experiment we use \t{funmv} (option 5) to compute the combination
$y(t)=\cos(t\sqrtA)b+t\,\sinc(t\sqrtA)z$.
Note that \t{trigmv} is inapplicable for this problem because
it requires an explicit computation of possibly dense $\sqrtA$.
The computation of a matrix square root is a challenging
problem itself and infeasible for large scale matrices.
Another difficulty is that \t{trigmv} cannot immediately yield
$x:=\sinc(t\sqrtA)b$, yet $x$ requires solving the system
$\sqrtA\,x=\sin(\sqrtA)b$, which could be dense or ill--conditioned.

Thus we invoke our \alg\ for the matrix $B = [b,z]$.
The combination above can be viewed as an exact solution of
the system \eqref{ode} with $g\equiv0$, $y(0)=b$ and $y'(0)=z$.
We compare the approximation of $y(t)$ using our \alg\ with that
obtained from the formula
\begin{equation}\label{blk.exp}
  \exp\left(t\left[
              \begin{array}{cc}
                0 & I \\
                -A & 0 \\
              \end{array}
             \right]\right)
             \left[
                     \begin{array}{c}
                       b \\
                       z \\
                      \end{array}\right]
                   = \left[
                   \begin{array}{c}
                       y(t) \\
                       y'(t) \\
                     \end{array}\right],
\end{equation}
which is a particular case of the expression given in
\cite[Prob.~4.1]{high:FM}; see also \cite[Eq.~(1.1)]{hika17}. We use the \Alg\ of Al-Mohy and Higham
\t{expmv} to evaluate the left hand side of \eqref{blk.exp}.
The approximation of $y(t)$ is obtained by reading off
the upper half of the resulting vector.
\begin{table}[t]  


\caption{\Experiment\ \ref{ex3}:
         \speed\ denotes time for method divided by
         time for \t{funmv}.}

\label{table.ex3}

\begin{center} \footnotesize

(a) Double precision\\
\begin{tabular}{cc|*{3}{c}|*{3}{c}|*{1}{c}}
\\\hline
&& \multicolumn{3}{c|}{\t{funmv} } &
   \multicolumn{3}{c|}{\t{expmv}} &
   \multicolumn{1}{c}{\t{expm}}\\
& $t$ & \speed  & \mv  & Error  & \speed & \mv   & Error & \speed\\
\hline
\mystrut{8pt}

\t{orani678}    &  100  &   1 &      920 &   3.2e-14  &   1.3 &     2046 &   3.2e-14&   7.1e2 \\
\t{bcspwr10}    &   10  &   1 &      190 &   4.5e-15  &   2.6 &      616 &   4.4e-15&   1.6e4 \\
\verb"gr_30_30" &    2  &   1 &       86 &   2.0e-15  &   1.5 &      180 &   1.7e-15&   3.3e2 \\
\t{triw}        &   10  &   1 &     1694 &   3.3e-14  &   3.8 &     4144 &   3.9e-14&   3.3e1
\end{tabular}
\vspace{0.3cm}
\\
(b) Single precision\\
\begin{tabular}{cc|*{3}{c}|*{3}{c}|*{1}{c}}
\\\hline
&& \multicolumn{3}{c|}{\t{funmv} } &
   \multicolumn{3}{c|}{\t{expmv}} &
   \multicolumn{1}{c}{\t{expm}}\\
& $t$ & \speed  & \mv  & Error  & \speed & \mv   & Error & \speed\\
\hline
\mystrut{8pt}
\t{orani678}    &  100  &   1 &      558 &   1.5e-9  &   1.5 &     1348 &   2.8e-8&   1.3e3 \\
\t{bcspwr10}    &   10  &   1 &      134 &   6.9e-11  &   3.0 &      496 &   2.3e-10&   2.3e4 \\
\verb"gr_30_30" &    2  &   1 &       58 &   3.2e-11  &   1.4 &       96 &   1.9e-10&   4.6e2 \\
\t{triw}        &   10  &   1 &      930 &   3.9e-11  &   3.6 &     2216 &   4.7e-9&   5.8e1
\end{tabular}
\vspace{0.3cm}
\\
(c) Half precision\\
\begin{tabular}{cc|*{3}{c}|*{3}{c}|*{1}{c}}
\\\hline
&& \multicolumn{3}{c|}{\t{funmv} } &
   \multicolumn{3}{c|}{\t{expmv}} &
   \multicolumn{1}{c}{\t{expm}}\\
& $t$ & \speed  & \mv  & Error  & \speed & \mv   & Error & \speed\\
\hline
\mystrut{8pt}
\t{orani678}    &  100  &   1 &      400 &   1.2e-4  &   1.6 &      992 &   1.8e-3&   1.8e3 \\
\t{bcspwr10}    &   10  &   1 &       80 &   2.4e-6  &   4.4 &      444 &   3.0e-5&   3.6e4 \\
\verb"gr_30_30" &    2  &   1 &       46 &   2.1e-7  &   1.1 &       60 &   1.8e-5&   5.2e2 \\
\t{triw}        &   10  &   1 &      650 &   6.9e-6  &   3.1 &     1370 &   2.5e-4&   8.0e1
\end{tabular}
\end{center}

\end{table}
For a reference solution we use the MATLAB function
\t{expm} to compute the left hand side of \eqref{blk.exp}.
We use the matrices and the vectors $b$ prescribed in Experiment \ref{ex2} except \t{poisson} due to memory
limitation because of the use of \t{expm}.
We take $z=[\sin 1, \sin 2,\cdots,\sin n]^T$ for all matrices.
For fairer comparison we multiply by two the number of matrix--vector
products \mv\ counted by the code \t{expmv} because the dimension of
the input matrices is $2n\times2n$. Table \ref{table.ex3} presents the results.
Obviously our \alg\ outperforms the alternative block version of
the problem in terms of CPU time and computational cost with slightly
better relative forward errors for single and half precisions.
Using the MATLAB function \t{profile} to analyze the execution time
for \t{funmv} and \t{expmv} in the experiment as a whole, the CPU time of \t{funmv}
represents around 22 percent of the CPU time of both functions.
\end{example}
\section{Concluding remarks}\label{sec5}
The \alg\ we developed here has direct applications to solving
second order systems of ODE's and their trigonometric numerical
schemes. A single invocation of \Alg\ \ref{alg.funmv} for inputs
$h$, $A$, and $B=[y_n,y'_n,\widehat{g}(y_n)]$ returns the six vectors
$\cos(h\sqrtA)y_n$, $\cos(h\sqrtA)y'_n$, $\cos(h\sqrtA)\widehat{g}(y_n)$,
$\sinc(h\sqrtA)y_n$, $\sinc(h\sqrtA)y'_n$, and $\sinc(h\sqrtA)\widehat{g}(y_n)$
that make up the vectors $y_{n+1}$ and $y'_{n+1}$ in the scheme
\eqref{yn} and \eqref{y'n}.
The evaluation of this scheme draws our attention back to the end of section
\ref{sec2}. Since the \alg\ has to be executed repeatedly for a fixed
matrix $A$ and different $B$ and perhaps different scalar $h$, it is
recommended to precompute the matrix $S_{pm}$ \eqref{spm} and provide it
as an external input to reduce the cost of the whole computation.

\Alg~\ref{alg.funmv} has several features. First, it computes the action of
the composition $f(t\sqrtA)B$ without explicitly computing $\sqrtA$.
Second, it returns results in finite number of steps that can be
predicted before executing the main phase of the \alg. Third, the \alg\ is
easy to implement and works for any matrix and the only external parameter
that control the computation is $\tol$. Fourth, the \alg\ spends most of
its work on multiplying $A$ by vectors. Thus it fully benefits from the
sparsity of $A$ and fast implementation of matrix multiplication.
Fifth, we can use \Alg~\ref{alg.funmv} (option 2) to compute the action
of the matrix exponential since $e^AB=\cosh(A)B+\sinh(A)B$.
Finally, though we derive the values of $\th_m$ in
\eqref{thm} for half, single, and double precisions, $\th_m$ can be
evaluated for any arbitrary precision. \Alg~\ref{alg.funmv} can be
extended to be a multiprecision \alg\ as in \cite{fahi17} since the
function $\rho_m$ \eqref{rhom} has an explicit expression that is easy
to be handled by optimization software.

All these features make our \alg\ attractive for black box use
in a wide range of applications.


\def\noopsort#1{}\def\l{\char32l}\renewcommand\v[1]{\accent20 #1}
  \let\^^_=\v\def\hbk{hardback}\def\pbk{paperback}

\end{document}